\documentclass{article}
\usepackage{amsmath,amssymb,amsthm,color}
\usepackage[dvips]{graphicx,epsfig}

\usepackage{ifthen}
\usepackage[latin1]{inputenc}
\usepackage{amsmath,amssymb,amsthm}
\usepackage[dvips]{graphicx,epsfig}
\usepackage{color} 

\def\AA{{\mathbb A}}

\def\NN{{\mathbb N}}

\def\PP{{\mathbb P}}

\def\RR{{\mathbb R}}

\def\ZZ{{\mathbb Z}}

\newcommand{\fd}{\ensuremath{\rightarrow}}

\newcommand{\inj}{\ensuremath{\hookrightarrow}}

\newcommand{\f}{\frac}
\newcommand{\findem}{\nolinebreak\vspace{\baselineskip} \hfill\rule{2mm}{2mm}\\}
\renewcommand{\phi}{\ensuremath{\varphi}}
\newcommand{\inc}{\ensuremath{\subset}}

\newcommand{\plp}{\PP^2}
\renewcommand{\phi}{\ensuremath{\varphi}}
\newcommand{\x}{\ensuremath{\times}}

\newcommand{\s}{Spec\ }

\newtheorem{nt}{Notation}
\newtheorem{conj}[nt]{Conjecture}
\newtheorem{prop}[nt]{Proposition}

\newtheorem{exercice}[nt]{Exercice} 
\newcounter{numeroquestion} 
\newcounter{numerosousquestion}

\newcommand{\sousquestion}{\ifthenelse{\value{numerosousquestion}=1}{}{\\}\textbf{\roman{numerosousquestion})} \addtocounter{numerosousquestion}{1}}

\newenvironment{listecompacte}
{\begin{list}
    {\ensuremath{\bullet}}
    {\setlength{\topsep}{2pt}
      \setlength{\itemsep}{1pt} \setlength{\parsep}{0pt}}
}
{\end{list}
}

\newtheorem{defi}[nt]{Definition}

\newtheorem{lm}[nt]{Lemma} 
 
\newtheorem{rem}[nt]{Remark} 
\newtheorem{thm}[nt]{Theorem}  
 
\newcommand{\demo}{\noindent \textit{Proof}}

\begin{document}
\sloppy
\bibliographystyle{plain}
\title{On the postulation of $s^d$ fat points in $\PP^d$ }
\date{}
\author{Laurent Evain (laurent.evain@univ-angers.fr)}
\maketitle
\noindent

\newcommand{\vecteur}[0]{\Delta=(\delta_1,\dots,\delta_d)}
\newcommand{\corvecteur}[0]{(\delta_1,\dots,\delta_d)}

\section*{Abstract}
In connection with his counter-example to the fourteenth problem of
Hilbert, Nagata formulated a conjecture concerning the postulation of 
$r$ fat points of the same multiplicity in $\plp$ and proved it when 
$r$ is a square. Iarrobino formulated a similar conjecture 
in $\PP^d$. We prove Iarrobino's conjecture when $r$ is a $d$-th power.
As a corollary, we obtain new counter-examples modeled on those by Nagata.

\section{Introduction}
\label{sec:Intro}
What is the dimension $l(d,\delta,\mu_1,\dots,\mu_r)$ 
of the sub-vector space of $k[X_0,\dots,X_d]$ containing the
homogeneous polynomials
of degree $\delta$ that vanish at general points $p_1,\dots,p_r \in \PP^d$
with order $\mu_1,\dots,\mu_r$ ? This question
remains open as soon as $d\geq 2$ and has numerous consequences
( see \cite{ciliberto-kouvidakis99}, \cite{iarrobino95:_waring_forms},
\cite{xu95:_ampleBundleOnSurfaces},
\cite{mcDuff-Polterovich94:symplecticPackingEtNagata},
\cite{nagata59:14emeProblemeDeHilbert}
for instance).

The question was raised by Nagata  in connection with his answer to the 
fourteenth problem of Hilbert
\cite{nagata59:14emeProblemeDeHilbert}. 
He gave an example of a linear
action on a finite dimensional vector space such that the algebra of  polynomial invariants
is not finitely generated. The key point in the proof,
which Nagata called the ``fundamental lemma'', is the equality 
$l(2,4m,m_1=m,\dots,m_{16}=m)=0$. 

When the dimension of the ambiant projective space is $d=2$
and the number of points is $r\leq 9$, 
the dimension  $l(d,\delta,\mu_1,\dots,\mu_r)$ is well known
\cite{nagata60:9points_imposent_conditions_indep}. 
As for the remaining cases $r>9$, 
Nagata formulated the following conjecture:
$$
l(2,\delta,\underbrace{\mu,\dots,\mu}_{r\ times})=l(2,\delta,\mu^r)
=0 \mbox{ if }\delta \leq \sqrt{r}\mu,
$$
and proved it when $r$ is a square. 
This conjecture is of particular interest since
it crystallizes the difficulties. Indeed, the expected dimension
$l(2,\delta,\mu^r)$ is
$max(0,v(2,\delta,\mu^r))$
where
\begin{displaymath}
v(2,\delta,\mu^r)=\f{(\delta+2).(\delta+1)}{2}-r.\f{\mu.(\mu+1)}{2}
\end{displaymath}
is the so-called virtual dimension. With any known method,
the hardest cases are the cases with $r$ fixed, $\mu>>r$ and 
the degree $\delta$ is such that the virtual dimension is zero. 
An immediate estimate shows that the critical $\delta$ for which 
the virtual dimension is zero is asymptotically equivalent to
$\sqrt{r}\mu$.
It follows that the hardest cases correspond to 
Nagata's conjecture. Nagata proved himself 
this conjecture when $r$ is a square.  

Leaving the two-dimensional case for the general case, there is still
a conjecture for the dimension $l(d,\delta,\mu_1,\dots,\mu_r)$, 
due to Iarrobino \cite{iarrobino97:conjectureGeneralisantNagata}
(see also \cite{lafaceUgaglia03:conjecturePostulationP3}).
Facing the critical cases too, 
he derived from his conjecture a generalisation of
Nagata's conjecture:
\begin{conj}
Let $(r,d)$ be a couple of integers with 
\begin{listecompacte}
  \item $d\geq 2$
  \item $r\geq max(d+5,2^d)$
  \item  $(r,d) \notin \{ (7,2),(8,2),(9,3)\}$.
\end{listecompacte}
If 
$\delta < \sqrt[d]r \mu$ then
$l(d,\delta,\mu^r)=0$. 
\end{conj}
In the 2-dimensional case however, this is not 
exactly Nagata's conjecture.
Indeed, Nagata's conjecture is very slightly stronger, 
since the condition on $\delta$ is $\delta \leq \sqrt r \mu$,
not $\delta < \sqrt r \mu$, and this difference 
turned out to be very important in the applications
(in Nagata's counter-example to the fourteenth problem
of Hilbert, or in \cite{ciliberto-kouvidakis99} for instance).
Replacing carelessly the strict inequality by a large 
inequality is not possible since the cases $(r,d)=(8,3)$ 
and $(r,d)=(9,2)$ would obviously contradict the statement. 
Nevertheless, excluding these cases, one can 
formulate the conjecture as follows:
\begin{conj}
Let $(r,d)$ be a couple of integers with 
\begin{listecompacte}
  \item $d\geq 2$
  \item $r\geq max(d+5,2^d)$
  \item $(r,d) \notin \{(7,2),(8,2),(9,2),(8,3),(9,3)\}$.
\end{listecompacte}
If
$\delta \leq \sqrt[d]r \mu$ then
$l(d,\delta,\mu^r)=0$.
\end{conj}
Let us call this conjecture the large critical conjecture in 
opposition to the conjecture by Iarrobino which we shall call the 
strict critical conjecture.

The goal of this paper is to prove that the large critical conjecture
holds when the number of points is a power with exponent the dimension of 
the ambiant projective space: 
\begin{thm} 
\label{thr:large}
  Let $k$ be an algebraically closed field of characteristic zero.
  Let $d\geq 2$ be an integer, $r$ be an integer such that $r=s^d$
  for some $s\geq 2$. Suppose moreover that $(r,d)\notin\{
  (4,2),(9,2),(8,3)\}$. Then:\\
  $l(d,\delta, \mu^{r})=0$ if $\delta \leq s\mu$.
\end{thm}
As a corollary, we obtain
new counter-examples to the fourteenth problem 
of Hilbert. Indeed, replacing
the fundamental lemma of Nagata with our theorem,
one can mimic step by step the construction of Nagata (with
a few minor and easy changes) to exhibit a new example. In concrete terms,
each couple $(s,d)$ of the theorem gives a new fundamental lemma and 
a new counter-example. The example associated with the
couple $(s,d)$ is an action of the affine group $G_a^{s^d-d-1}$
on a vector space of dimension $2s^d$:
\begin{thm}\label{nouveauContreExemple} 
Let $a_{ij}$ ($i=0\dots d,\ j=1,\dots,s^d$) be the coordinates of
$s^d$ generic points of $\PP^d$. Let $V$ be the vector space 
of dimension $s^d$ and $V^*\subset V$ be the set of vectors orthogonal 
to the $d+1$ vectors $(a_{i1},\dots,a_{is^d})$. Let $G$ be the set of 
linear transformations $\sigma$ of $\s
k[x_1,\dots,x_{s^d},t_1,\dots,t_{s^d}]$ such that  
\begin{listecompacte}
  \item $\sigma(t_i)=t_i$
  \item $\sigma(x_i)=x_i+b_i t_i$
\end{listecompacte}
for some $(b_1,\dots,b_{s^d})\in V^*$.
Then the algebra of elements of
$k[x_1,\dots,x_{s^d},t_1,\dots,t_{s^d}]$
invariant under $G$ is not finitely generated. 
\end{thm}
As mentioned, the proof of theorem \ref{nouveauContreExemple} is
a straightforward generalisation
of Nagata's proof \cite{nagata59:14emeProblemeDeHilbert} and we 
refer to this paper for it. 

Our method to prove theorem \ref{thr:large} is an induction 
on the dimension of the ambiant projective space.
The formulation of the theorem does not suggest 
such an induction; however, using the notion of collision of fat
points, we transform the statement of the theorem into a combinatorial 
statement and we perform the induction on the combinatorial
statement (see remark \ref{rem:induction on dimension}).

\begin{rem} \label{remarque sur le cas 2^d}
It seems that theorem \ref{thr:large} 
leaves the cases $(r,d)=(4,2)$, $(r,d)=~(9,2)$ 
and $(r,d)=(8,3)$ untreated. 
However, these  cases are completly understood.
Indeed, by 
\cite{nagata60:9points_imposent_conditions_indep} for $(r,d)=(4,2)$
and $(9,2)$,
and 
by proposition \ref{prop: le cas r=2^d} for $(r,d)=(8,3)$,
we have $l(d,\delta,\mu^r)=max(0,\binom
  {\delta+d}{d}-r.\binom{d+\mu-1}{d})$.
\end{rem}
If the characteristic of the base field is arbitrary,
we can forget the parts of the proof which 
use the hypothesis on the characteristic and we still have
the strict critical
conjecture:
\begin{thm} \label{thr:strict}
Let $d\geq 2$ be an integer and let $r$ be a $d^{th}$-power.
If $\delta < \sqrt[d]r \mu$ then
$l(d,\delta,\mu^r)=0$.
\end{thm}

\section{Stratifications on the Hilbert scheme}
\label{sec: def des stratifications}
In this section, we explain the strategy of the proof: we
define locally closed subschemes 
$C(E_1,\dots,E_i)$ of the Hilbert 
scheme $Hilb(\PP^d)$ 
and we reduce the proof to an incidence between these
subschemes. 

\subsubsection*{Monomial subschemes}
\label{sec:monomial-subschemes}
A staircase $E$ in $\NN^d$ is a subset whose 
complementary $\NN^d-E$ verifies 
$$(\NN^d-E)+\NN^d\inc \NN^d-E.$$
A staircase $E$ being fixed, let $I^E \inc k[[x_1,\dots,x_d]]$
(resp. $I^E \inc k[x_1,\dots,x_d]$)
be the ideal whose elements are  
the series (resp. the polynomials)
$$
\sum c_{\alpha_1\alpha_2\dots\alpha_d}
x_1^{\alpha_1}x_2^{\alpha_2}\dots
x_d^{\alpha_d}=\sum c_{\underline \alpha}{\underline x}^{\underline
  \alpha}
$$
verifying $c_{\underline\alpha}=0$ if $\underline \alpha \in E$. A 
zero-dimensional subscheme $Z$
of $\PP^d$ supported by a point $q$ 
is said to be monomial with staircase 
$E$ if it is defined by the ideal $I^E$ in 
a suitable formal neighborhood $\s k[[x_1,\dots,x_d]]\inj \PP^d$ of
$q$.
\\
A fat point of multiplicity $m$ 
is by definition a monomial subscheme defined by the 
regular staircase $R_m$: 
$$R_m:=\{(\alpha_1,\dots,\alpha_d) \ s.t.\ \alpha_1+\dots+\alpha_d<m\}.$$

\subsubsection*{Subschemes of $Hilb(\PP^d)$}
\label{sec:subschemes-hilbppd}

If $E_1,\dots,E_i$ are finite staircases in $\NN^d$, we denote by 
$C(E_1,\dots,E_i)$ the reduced subscheme of $Hilb\ \PP^d$ 
whose points parametrize the subschemes $Z$ of 
$\PP^d$ which are the disjoint 
union of $i$ distinct monomial subschemes
with staircases $E_1,\dots,E_i$. In symbols
$Z=\coprod Z_j$, where $Z_j$ is monomial with staircase $E_j$. 
It is known by \cite{hirscho85:methodeHoraceManuscripta} 
that $C(E_1,\dots,E_i)\inc Hilb \PP^d $ is 
a locally closed irreducible subscheme. In particular it has 
a generic point $G$, which parametrizes a subscheme $Z_G$
whose ideal is denoted by $I_{Z_G}$. 
We denote by $l(d,\delta,E_1,\dots,E_i)=h^0(I_{Z_G}(\delta))$
the number of independant hypersurfaces of degree $\delta$
in $\PP^d$ containing $Z_G$.

\subsubsection*{Iarrobino's conjecture and incidence between strata}
\label{sec:reform-iarr-conj}

The theorem we want to prove can obviously be reformulated
as:
\begin{thm}
Let $r=s^d$ and $\delta \leq s\mu$.
Then:\\
$l(d,\delta,
\underbrace{R_{\mu},\dots,R_{\mu}}_{r\ times})=0$
if $(s,d) \notin \{(1,d),(2,d),(3,2)\}$ and if the 
characteristic of the base 
field is zero.
\end{thm}
The following proposition reduces the proof of the theorem to 
the computation of the closure of $C(R_{\mu},\dots,R_{\mu})$.

\begin{prop} \label{reduction_to_the_computation_of_closure}
  Let $E_1,\dots,E_i \inc \NN^d$ be staircases. Suppose that 
there exists a staircase $F$ with $F\supset R_{\delta+1}$
and $C(F)\inc \overline {C(E_1,\dots,E_i)}$,
then $l(d,\delta,E_1,\dots,E_i)=~0$. 
\end{prop}
\demo:
by semi-continuity of the cohomology
$l(d,\delta,E_1,\dots,E_i)\leq l(d,\delta,F)$
and $l(d,\delta,F)\leq l(d,\delta,R_{\delta+1})$ since $F\supset R_{\delta+1}$.
Since obviously $l(d,\delta,R_{\delta+1})=0$, the vanishing 
of $l(d,\delta,E_1,\dots,E_i)$ follows from the last
two inequalities.
\findem

\section{Elementary Incidences}
\label{sec: description des incidences elementaires}
The previous section explained that the theorems would follow
from incidences between the various subschemes $C(E_1,\dots,E_j)$. 
The goal of this section is to exhibit such incidences. 

Let $E \inc \NN^d$ be a finite staircase and $i\in \{1,\dots,d\}$ be an integer. 
There exists a unique ``height'' function 
$$
h_{E,i}:\NN^{d-1}\fd \NN
$$
such that 
$$(a_1,\dots,a_d)\in E \Leftrightarrow
a_i<h_{E,i}(a_1,\dots,a_{i-1},a_{i+1},\dots,a_d)
$$
Conversely, a function $h$ is the height function
of some staircase if and only if $h(a+b)\leq h(a)$ 
for any $(a,b)\in \NN^{d-1}\x \NN^{d-1}$.  
If $E_1,\dots,E_j$ are staircases, the sum of 
$E_1,\dots,E_j$ along the $i^{th}$ coordinate is 
the staircase $S_i(E_1,\dots,E_j)$ characterized 
by its height function 
$$
h_{S_i(E_1,\dots,E_j),i}=\sum_{k=1}^j h_{E_k,i}.
$$

\begin{prop} \label{prop:incidence par collision de front}
  Let $E_1,\dots,E_j$ be staircases and $k \in \{1,\dots,j\}$.
  Then $\overline {C(E_1,\dots,E_j)} \supset C(S_i(E_1,\dots
  E_k),E_{k+1},\dots,E_j)$. 
\end{prop}
\demo: this is a straightforward generalisation of 
\cite{hirscho85:methodeHoraceManuscripta}, proposition 5.1.2.
\findem

Let $(a_1,\dots,a_d)\in (\NN^*)^d$ and let $E$ be a staircase. 
We denote by $(a_1,\dots,a_d).E$ the staircase ``obtained from 
$E$'' by the linear map 
\begin{displaymath}
(x_1,\dots,x_d)\mapsto (a_1x_1,\dots,a_dx_d).
\end{displaymath}
Concretely, this is the smallest staircase 
satisfying the relation:
$$
(m_1,\dots,m_d) \in E 
\Rightarrow
(a_1(m_1+1)-1,\dots,a_d(m_d+1)-1)\in (a_1,\dots,a_d).E
$$
This is a staircase of cardinal $a_1.a_2\dots a_d.\#E$.
We denote by $a.E$ the staircase $(a,a,\dots,a).E$.

\begin{prop} \label{prop: incidence par affinite}
  Let $E,E_1,\dots,E_j$ be staircases. Then:\\
  $\overline {C(\underbrace{E,\dots,E}_{\prod a_i\ times},E_1,\dots,E_j)}
  \supset C((a_1,\dots,a_d).E,E_1,\dots,E_j)$.
\end{prop}
\demo: by induction on the number of $a_i$'s which are not equal to
one.  If all the $a_i$'s  but one are equal to one, the statement
follows from the previous proposition since
$$(1,\dots,1,a_i,1,\dots,1).E=S_i(\underbrace{E,\dots,E}_{a_i \mbox{
    times }}).
$$
For the general case, one can suppose by symmetry that $a_1\neq 1$. 
Applying several times -namely $a_2 a_3 \dots a_d$ times- this first step, we get
$$
\overline {C(\underbrace{E,\dots,E}_{\prod a_i\ times},E_1,\dots,E_j)}
\supset 
C(\underbrace{(a_1,1,\dots,1).E,\dots,(a_1,1,\dots,1)..E}_{a_2\dots a_d
  \mbox { times }},E_1,\dots,E_j)
$$
and, by induction,  
$$
\overline {C(\underbrace{(a_1,1,\dots,1).E,\dots,(a_1,1,\dots,1).E}_{a_2\dots a_d
  \mbox { times }},E_1,\dots,E_j)} 
$$
contains
$$
C((1,a_2,\dots,a_d).(a_1,1,\dots,1).E,E_1,\dots,E_j)=C((a_1,\dots,a_d).E,E_1,\dots,E_j).
$$
The expected inclusion follows immediatly.
\findem

In particular, when $a_1=a_2=\dots=a_d=s$, we get: 
\begin{prop} \label{corollaire: incidence par homotetie}
  Let $E,E_1,\dots,E_j$ be staircases. Then:\\
  $\overline {C(\underbrace{E,\dots,E}_{s^d\ times},E_1,\dots,E_j)}
  \supset C(s.E,E_1,\dots,E_j)$.
\end{prop}

\begin{defi}
  \label{def:1}Let $\vecteur$ be a primitive vector in $\ZZ^d$ such that there exist $i,j$ 
  satisfying $\delta_i\delta_j<0$. Let $E\inc \NN^d$ be a
  subset. We denote by 
  $\Delta(E) \inc \NN^d$ the unique subset verifying the following
  two conditions:
  \begin{listecompacte}
  \item for any line $L$ in $\RR^d$ with direction $\Delta$, the sets 
    $E \cap L$ and $\Delta(E)\cap L$ are equipotent
  \item $\forall i\in \NN$, $\forall (n,p) \in (\NN^d)^2$, 
$n\in \Delta(E)$ and $p=n+i\Delta \Rightarrow p \in \Delta(E)$
  \end{listecompacte}
\end{defi}
To be more explicit, the set $L\cap \NN^d$ is finite by hypothesis on 
$\Delta$. If $m_1<m_2<\dots<m_j$ are its elements, 
ordered by the relation
$$
(<)\hspace{1cm}m_{i_1}<m_{i_2} \Leftrightarrow \exists i\in \NN,\ m_{i_1}=m_{i_2}+i\Delta,
$$
then $\Delta(E)\cap L=\{m_1,\dots,m_k\}$, where $k=\#(E \cap L)$.

\begin{prop}     \label{prop:incidence pour escalier unique}
  Let $\vecteur \in \ZZ^d$ be a vector such that
  \begin{listecompacte}
    \item $\exists i,\ \delta_i=1$,
    \item $\forall k, \ k\neq i \Rightarrow \delta_k\leq 0$,
    \item $\exists j\neq i, \delta_j \neq 0$.
  \end{listecompacte}
  Then for every staircase $E$, $\Delta(E)$ is a staircase.
  Moreover, we
  have in characteristic zero
  the incidence:\\  $\overline{C(E)} \supset C(\Delta(E))$
\end{prop}
\demo: suppose by symmetry that $\delta_1 =1$. Let 
\begin{eqnarray*}
 \Phi:k[x_1,\dots,x_d] &\fd& k[x_1,\dots,x_d][t,\f{1}{t}]\\
x_1 & \mapsto & tx_1 + x_2^{-\delta_2}x_3^{-\delta_3}\dots x_d^{-\delta_d}\\
x_i & \mapsto & x_i \mbox{ if $i\neq 1$}. 
\end{eqnarray*}
The ideal 
$$I(t)=k[x_1,\dots,x_d][t,\f{1}{t}]\Phi(I^E)$$
defines a subscheme $$F\subset (\AA^1-\{0\})\x \AA^d$$ 
whose fiber over each $t \in \AA^1-\{0\}$ is a monomial 
subscheme with staircase $E$. In particular, $F$ is flat over 
$\AA^1-\{0\}$. The closure $\overline F \inc \AA^1 \x \AA^d$
is defined by the ideal $J(t)=I(t)\cap k[x_1,\dots,x_d,t]$ and it is 
flat over $\AA^1$. 

We want to prove the equality $J(0)=I^{\Delta(E)}$, using 
a natural graduation. 

Let $\phi_1,\dots,\phi_{d-1}:\ZZ^d\fd \ZZ$ 
be independant linear forms which vanish on $\Delta$. Consider the
multi-graduation $D$ defined by:
\begin{eqnarray*}
  D:\mathrm{Monomials\ of\ }k[x_1,\dots,x_d] &\fd& \ZZ^{d-1}\\
\underline{x}^{\underline{\alpha}} &\mapsto& 
(\phi_1(\underline{\alpha}),\dots,\phi_{d-1}(\underline{\alpha}))
\end{eqnarray*}
The conditions on $\Delta$ imply that, for all $\underline z=(z_1,\dots,z_{d-1})\in
\ZZ^{d-1}$, the sub-vector space
$k[x_1,\dots,x_{d}]_{\underline z}\inc k[x_1,\dots,x_{d}]$
containing the elements of degree $\underline z$ has finite
dimension. 
Note that $J(t)$ is a graded ideal ie. 
$$
J(t)=\oplus_{\underline z \in \ZZ^{d-1}} J_{\underline z}(t) 
$$
where 
$$J_{\underline z}(t)=J(t)\cap k[x_1,\dots,x_d]_{\underline z}[t].$$
In particular,
to compute $J(0)=\lim_{t\fd 0}J(t)$,
it suffices to compute the limit of its graded parts
in 
the grassmannians $G(l,k[x_1,\dots,x_d]_{\underline z})$
, where $l=\dim J_{\underline z}(t), t\neq 0$.
Let $m_1<\dots<m_k$ be the monomials of $k[x_1,\dots,x_d]_{\underline
  z}$,
where the order is given by the relation $(<)$ above. 
Let us admit temporarily
the inclusion
$$ (*) \  m_{k-l+1},m_{k-l+2},\dots,m_{k} \in J_{\underline z}(0). 
$$
Then $J_ {\underline z}(0)$ is the vector space generated by 
$ m_{k-l+1},m_{k-l+2},\dots,m_{k} $ for dimensional reasons
and $J(0)=I^{\Delta(E)}$ since these two graded ideals have 
the same graded parts.
In particular $J(0)$ is an ideal generated by monomials and the set 
$\Delta(E)$ of monomials which are not in 
$J(0)$ is a staircase. Moreover, replacing the coordinates 
$x_1,\dots,x_d$ of $\AA^d$ by any local system of coordinates, 
one shows by the same computation that any closed point of 
$C(\Delta(E),E_1,\dots,E_j)$ is a limit of points which are in 
$\overline{C(E,E_1,\dots,E_j)}$. This gives the incidence 
between the strata.  

It remains to show $(*)$.
Let $n_1=x^{{\underline \alpha}(1)},\dots,n_l=x^{{\underline \alpha}(l)}$ 
be the monomials of $I^E\cap k[x_1,\dots,x_d]_{\underline z}$, where 
$\underline \alpha(i)=(\alpha_1(i),\dots,\alpha_d(i))$.
The ideal $I(t)$
contains the monomials
\begin{eqnarray*}
\Phi(n_i)&=& (tx_1+x_2^{-\delta_2}x_3^{-\delta_3}\dots x_d^{-\delta_d}
)^{\alpha_1(i)}
x_2^{\alpha_2(i)}\dots
x_d^{\alpha_d(i)}.
\end{eqnarray*}
Since 
the degree of $m_i$ in $x_1$ is $k-i$,
this equality can be rewritten as:
\begin{eqnarray*}
\Phi(n_i)
=\sum_{j=0}^{\alpha_1(i)} \binom{\alpha_1(i)}{j} t^{j} m_{k-j}
&=& \sum_{j=0}^{k-1} \binom{\alpha_1(i)}{j} t^{j} m_{k-j}
\end{eqnarray*}
with the usual convention $\binom{\alpha_1(i)}{j}=0$ if
$j>\alpha_1(i)$. If $N$ and $M$ are the column matrices whose 
entries are respectively $\Phi(n_i)$, $i\in \{1,\dots,l\}$, and 
$t^{j}m_{k-j}$, $j~\in~\{0,\dots,k-1\}$, if $P$ is the
matrix whose coefficient $P_{ij}$ is $\binom {\alpha_1(i)}{j}$, 
the above equality writes down $N=PM$. Take the first $l$
columns of $P$ to get a square matrix 
$$Q=
\left(
  \begin{array}{ccccc}
    1 & \alpha_1(1) & \binom{\alpha_1(1)}{2} & \dots &
    \binom{\alpha_1(1)}{l-1}\\
    \dots & \dots \\
    1 & \alpha_1(l) & \binom{\alpha_1(l)}{2} & \dots &
    \binom{\alpha_1(l)}{l-1}
  \end{array}
\right).
$$
Since the coefficients in the third column are polynomials of degree 2 
in $\alpha_1$, one can replace the third column by a linear
combination of the first three columns
so that the $i^{th}$ element in the third column 
becomes $\alpha_1(i)^2$. Similarly, after suitable operations on the
columns, the $i^{th}$ element in the fourth, fifth column \ldots 
becomes $\alpha_1(i)^3,\alpha_1(i)^4,\dots $.
The resulting matrix is a Van Der Monde matrix in the $\alpha_1(i)$'s.
In characteristic zero, its determinant is
not zero since the $\alpha_1(i)$'s are distinct. In particular $Q$ is 
invertible. 

The ideal $I(t)$ contains the elements which are the coefficients 
of the matrix $Q^{-1}N=Q^{-1}PM$. Using that the identity is a submatrix 
of $Q^{-1}P$ by construction, the $i^{th}$ 
element in this column matrix is $c_i(t)=t^{i-1}m_{k-i+1}+R$ where $R$ is a
polynomial dividible by $t^{i}$. Thus, 
$\f{c_i(t)}{t^{i-1}} \in J(t)$ and, 
as expected, $J(0)$ contains $\f
{c_i(t)}{t^{i-1}}(0)=m_{k-i+1}$ for $i\in \{1,\dots,l\}$.
\findem
\\
If we have a finite set of monomial subschemes, we can specialize the
first one and leave the remaining subschemes unchanged. Thus, we get 
as a corollary of the previous proposition:
\begin{prop} \label{prop:incidence par specialisation le long de droites}
  Let $\vecteur \in \ZZ^d$ be a vector such that
  \begin{listecompacte}
    \item $\exists i,\ \delta_i=1$,
    \item $\forall k, \ k\neq i \Rightarrow \delta_k\leq 0$,
    \item $\exists j\neq i, \delta_j \neq 0$.
  \end{listecompacte}
  Then for every set of staircases $E,E_1,\dots,E_j$, we
  have in characteristic zero
  the incidence:\\  $\overline{C(E,E_1,\dots,E_j)} \supset C(\Delta(E),E_1,\dots,E_j)$
\end{prop}

\subsection{Combinatorial properties of $\Delta$ }
\label{sec:comb-prop-delta}
We give here some combinatorial properties of the map $E \mapsto \Delta(E)$
that we will use later on. 

\begin{lm}
Let $E$ and $F$ be two subsets of $\NN^d$ and $\vecteur \in \ZZ^d$ be a direction 
satisfying the properties of the preceding proposition. Suppose
that for every line $L$ with direction $\Delta$, we have the
inequality on cardinals:
$$ \#\{E\cap L\} \geq \#\{F \cap L\}$$
then $\Delta(E) \supset \Delta(F)$.
\end{lm}
\demo: we must show for everery line $L$ the inclusion
$\Delta(E)\cap L\; \supset \Delta(F)\cap L$.
This is obvious since, using
the $m_i$'s introduced after definition \ref{def:1},
$\Delta(E)\cap L=\{m_1,\dots,m_{\#\{E\cap L\}}
 \}$ and 
$\Delta(F)\cap L=\{m_1,\dots,m_{\#\{F \cap L\}}\}$. 
\findem

Applying this lemma to the following $E$ and to $F=R_{\mu}$,
noticing that $\Delta(R_{\mu})=R_\mu$, we get: 
\begin{lm}
\label{on bouche le dernier trou}
  Let $R_\mu$ be a regular staircase, $m\in R_{\mu}$, $P \inc \NN^d$ 
a subset such that $P\cap R_\mu=\emptyset$ and $E=R_{\mu}\cup P - \{m\}$. 
If there exists $i \in \ZZ$ such that $m+i \Delta \in P$, then 
$\Delta (E) \supset R_\mu$.
\end{lm}

\begin{lm}
\label{lm:lemme cle}
  Let $(s,d)$ be a couple of integers with $d\geq 2$,$s\geq
  2$, and $(s,d)\notin \{(2,2),(2,3),(3,2)\}$. Then
  there exists ($\Delta_d,\dots,\Delta_1) \in (\ZZ^d)^d$ such that
  \begin{listecompacte}
  \item $\forall i$, $\Delta_i$ verifies the conditions of proposition
    \ref{prop:incidence pour escalier unique}
  \item $\forall \mu>0$, $\Delta_d(\Delta_{d-1}(\dots(\Delta_1(s.R_{\mu}))))\supset
    R_{s\mu+1}$. 
  \end{listecompacte}
\end{lm}
\begin{rem}
  More precisely, it will follow from the proof that the choice of the
  $\Delta_i$ depend on $s$ in the following way. 
  \begin{listecompacte}
    \item $s>3$: $\Delta_1=(0,\dots,0,1,-s+1)$,
    $\Delta_2=(0,\dots,0,-s+2,1)$,
    $\Delta_i=(0,\dots,0,1,-1,-1,0,\dots,0)$ for $i\geq 3$, where the
    $1$ is on the position of index $1+d-i$. 
  \item $s=3$: $\Delta_1=(0,\dots,0,1,-2,0)$,
    $\Delta_2=(0,\dots,0,-3,0,1)$, $\Delta_3=(0,\dots,0,0,1,-2)$,
    $\Delta_i=(0,\dots,0,1,-1,-1,0,\dots,0)$ for $i\geq 4$.
  \item $s=2$: $\Delta_1=(0,\dots,0,1,-1,-1,-1)$,
    $\Delta_2=(0,\dots,0,-1,1,-1,0)$,
    $\Delta_3=(0,\dots,0,-1,0,1,-1)$,
     $\Delta_4=(0,\dots,0,-1,-1,0,1)$,
    $\Delta_i=(0,\dots,0,1,-1,-1,0,\dots,0)$ for $i\geq 5$.
  \end{listecompacte}

\end{rem}
\demo. We proceed by induction on $d$. 
Considering the couples $(s,d)$ involved in the proposition, 
we have to initialize the induction with the cases $(s>3,d=2)$,
$(s=3,d=3)$ and $(s=2,d=4)$. \\
\underline{Initial cases.}
If $d=2,s>3$, then one can take $\Delta_1=(1,-s+1)$ and
$\Delta_2=(-s+2,1)$. When $s=3,d=3$,
we must find $\Delta_1,\Delta_2,\Delta_3$ 
such that 
$$
\Delta_3(\Delta_2(\Delta_1(3.R_{\mu}))) \supset R_{3\mu+1}.
$$
The $\f{(\mu+1)(\mu+2)}{2}$ elements of the difference
$$
R_{3\mu+1}-3.R_{\mu}=\{(3x,3y,3z),\ x+y+z=\mu\}
$$
are shown in the following 
figure 

  \begin{figure}[h] 
     \begin{center}
        \input{specialisationEnDim3.pstex_t}
     \end{center} 
  \end{figure} 
with $\mu=2$. Taking $\Delta_1=(1,-2,0)$, we have:
$$
R_{3\mu+1}-\Delta_1(3.R_{\mu})=\{(0,3y,3z),\ y+z=\mu\}
$$
Finally, taking $\Delta_2=(-3,0,1)$ and $\Delta_3=(0,1,-2)$,
$$
\Delta_3(\Delta_2(\Delta_1(3.R_{\mu}))) \supset R_{3\mu+1},
$$
as expected.\\
Consider now the last initial case $(s=2,d=4)$. 
By definition, 
\begin{displaymath}
  2R_\mu=\{(x,y,z,t)\ s.t.\ [\frac{x}{2}] +[\frac{y}{2}] +[\frac{z}{2}]
  +[\frac{t}{2}] < \mu\}
\end{displaymath}
where $[\ ]$ stands for the integral part. 
If $P\inc \NN^4$ is a subset, we denote by $P^i$ the subset of $P$
containing the elements $(x,y,z,t)$ such that $i$ elements among
$(x,y,z,t)$ are odd and we put
$S_m=R_{m+1}-R_m$. With these notations, easy considerations on 
the parities of $(x,y,z,t)$ 
give the equality:
\begin{displaymath}
  2R_\mu=R_{2\mu}\coprod (S_{2\mu}\setminus S_{2\mu}^0) \coprod S_{2\mu+1}^3
  \coprod S_{2\mu+2}^4.
\end{displaymath}
To compute $\Delta(2R_\mu)$, we note that we can define $\Delta$ on
subsets in such a way that if $E=\coprod E_i$ is a disjoint union,
then $\Delta(E)=\coprod \Delta(E_i)$. Indeed, by construction
of the map $\Delta$, if $L$
is a line in $\RR^d$ with direction $\Delta$, then $E\cap L$ and
$\Delta(E)\cap L$ are two totally ordered sets of the same finite
cardinality, hence there is a unique increasing one-to-one
correspondance between $E
\cap L$ and $\Delta(E)\cap L$. If $e\in E$ and $L$ is the line with
direction $\Delta$ passing through $e$, $\Delta(e)$ is the
image of $e$ through this correspondance. We let
$\Delta(E_i)=\cup_{e_i\in E_i} \Delta(e_i)$.
\\
Let $\Delta_1=(1,-1,-1,-1)$. To compute the image $\Delta_1(e)$ of an
element $e$, we make the following observation. 
If $E\inc \NN^4$ can be written as a disjoint union
\begin{displaymath}
E=R_{j+1}\coprod E_{j+1} \coprod  E_{j+2} \coprod\dots E_{j+k},\mbox{ with\ } E_l\subset S_l
\end{displaymath}
and if $\Delta_1=(\Delta_x,\Delta_y,\Delta_z,\Delta_t)$ satisfies 
$-2(\Delta_x+\Delta_y+\Delta_z+\Delta_t)\geq k$, then
\begin{eqnarray*}
  \Delta_1(e)&=&e+\Delta_1  \mbox{\  \ if }  e+\Delta_1 \in \NN^4\setminus E\\
\Delta_1(e)&=&e \mbox{\ \ otherwise }. 
\end{eqnarray*}
This observation leads to the equality 
\begin{displaymath}
\Delta_1(2R_\mu)= R_{2\mu}\coprod (S_{2\mu}\setminus S_{2\mu}^0) \coprod S_{2\mu+1}^3
  \coprod \Delta_1(S_{2\mu+2}^4).
\end{displaymath}
If $P\inc \NN^4$, we define $P(1,*,\neq 0,e)\inc P$ to be the subset
containing the elements $(x,y,z,t)$ with $x=1$, $y$ any number, 
$z\neq 0$ and $t$
even. There are obvious generalisations of this notation. 
With this notation, we
have:
\begin{eqnarray*}
\Delta_1(2R_\mu)&=& R_{2\mu}\coprod (S_{2\mu}\setminus S_{2\mu}^0) \coprod S_{2\mu+1}^3
\coprod    S_{2\mu}^0(\neq 0,*,*,*)\\
&=&  R_{2\mu}\coprod (S_{2\mu}\setminus S_{2\mu}^0(0,*,*,*)) \coprod S_{2\mu+1}^3\\
&=& R_{2\mu} \coprod (S_{2\mu}\setminus S_{2\mu}^0(0,*,*,*)) \coprod
S_{2\mu+1}^3(1,*,*,e)  \\&& \coprod
S_{2\mu+1}^3(1,*,e,*) \coprod  S_{2\mu+1}^3(1,e,*,*) \coprod
S_{2\mu+1}^3(\neq 1,*,*,*).
\end{eqnarray*}
Let $\Delta_2=(-1,1,-1,0)$,
    $\Delta_3=(-1,0,1,-1)$,
     $\Delta_4=(-1,-1,0,1)$.
Then,
\begin{eqnarray*}
 \Delta_2 \circ \Delta_1(2R_\mu)&=& R_{2\mu} \coprod (S_{2\mu}\setminus S_{2\mu}^0(0,*,*,*)) \coprod
\Delta_2(S_{2\mu+1}^3(1,*,*,e))  \\&& \coprod
S_{2\mu+1}^3(1,*,e,*) \coprod  S_{2\mu+1}^3(1,e,*,*) \coprod
S_{2\mu+1}^3(\neq 1,*,*,*)
\\
&=& R_{2\mu} \coprod (S_{2\mu}\setminus S_{2\mu}^0(0,*,*,*)) \coprod
S_{2\mu}^0(0,\neq 0,*,*) \\&& \coprod
S_{2\mu+1}^3(1,*,e,*) \coprod  S_{2\mu+1}^3(1,e,*,*) \coprod
S_{2\mu+1}^3(\neq 1,*,*,*)
\\
&=& R_{2\mu} \coprod (S_{2\mu}\setminus S_{2\mu}^0(0,0,*,*)) \\&& \coprod
S_{2\mu+1}^3(1,*,e,*) \coprod  S_{2\mu+1}^3(1,e,*,*) \coprod
S_{2\mu+1}^3(\neq 1,*,*,*)
\end{eqnarray*}
\begin{eqnarray*}
 \Delta_3 \circ \Delta_2 \circ \Delta_1(2R_\mu)&=&  R_{2\mu} 
\coprod (S_{2\mu}\setminus S_{2\mu}^0(0,0,*,*)) \\&& \coprod
S_{2\mu+1}^3(1,*,e,*) \coprod  \Delta_3(S_{2\mu+1}^3(1,e,*,*))  \\&&\coprod
S_{2\mu+1}^3(\neq 1,*,*,*)
\\
&\supset &  R_{2\mu} 
\coprod (S_{2\mu}\setminus S_{2\mu}^0(0,0,0,*)) \\&& \coprod
S_{2\mu+1}^3(1,*,e,*) \coprod
S_{2\mu+1}^3(\neq 1,*,*,*)
\\
\Delta_4 \circ \Delta_3 \circ \Delta_2 \circ \Delta_1(2R_\mu)&\supset&  
 R_{2\mu} 
\coprod (S_{2\mu}\setminus (0,0,0,2\mu)) \\&& \coprod
\Delta_4(S_{2\mu+1}^3(1,*,e,*)) \coprod
S_{2\mu+1}^3(\neq 1,*,*,*)\\
&\supset &  R_{2\mu} 
\coprod S_{2\mu} \coprod
S_{2\mu+1}^3(\neq 1,*,*,*)
\end{eqnarray*}
We have obtained the required inclusion 
 $ \Delta_4 \circ \Delta_3 \circ\Delta_2 \circ \Delta_1(2R_\mu) \supset R_{2\mu+1}$.
\medskip 
\\
\underline{Step from $d-1$ to $d$}. 
As we will proceed by induction on the dimension $d$, 
we precise our notations and 
we denote by $R_i(d)$ the regular staircase $R_i$ in $\NN^d$. 
Let $T_i$ be the ``$i^{th}$ slice'' of
$s.R_{\mu}(d)$, i.e.
$$T_i:=\{m\in \NN^{d-1}\ s.t.\ (i,m)\in  s.R_{\mu}(d)\}.$$
Then $T_i=s.R_{{\nu(i)}}(d-1)$  with ${\nu(i)}=max(0,\mu-[\f{i}{s}])$.
When $i<s\mu$, $\nu(i)>0$ and we can apply the induction to  $T_i$.
Using moreover that $s\nu(i)\geq s\mu-i$, we get
elements $\gamma_1,\dots, \gamma_{d-1}
\in \NN^{d-1}$ such that, 
$$
T'_i=\gamma_{d-1}(\dots(\
\gamma_{1}(T_i))) \supset R_{s{\nu(i)}+1}(d-1)\supset R_{s\mu+1-i}(d-1).
$$
Let $\Delta_i=(0,\gamma_i) \in \NN^d$.
The  $i^{th}$ slice of the staircase 
\begin{displaymath}
F=\Delta_{d-1}(\dots(\Delta_1(s.R_{\mu}(d))))
\end{displaymath}
is $T'_i$. Summing up, for $i<s\mu$, the $i^{th}$ slice of $F$ strictly
contains the $i^{th}$ slice $R_{s\mu+1-i}(d-1)$ of $R_{s\mu+1}(d)$.
In particular, $F$
contains all the $d$-tuples whose sum is $s\mu$ except
$(s\mu,0,0,\dots,0)$.

It remains to find $\Delta_d$ such that $\Delta_d(F) \supset
R_{s\mu+1}(d)$ by an application of lemma
\ref{on bouche le dernier trou}.\\
Note that 
$$
T_{s\mu-1}=s.R_1(d-1)\supset K= R_2(d-1)\cup (1,1,0,\dots,0).
$$
It follows that 
$$
T'_{s\mu-1}\supset \gamma_{d-1}(\dots(\gamma_1(K)))=K
$$
and that the 
element $z=(s\mu-1,1,1,0,\dots,0)$
is in $F$. Let $\Delta_d=(1,-1,-1,0,\dots,0)$.
Applying lemma
\ref{on bouche le dernier trou}
with $m=(s\mu,0,0,\dots,0)$,
$E=F$, $P=F-R_{s\mu+1}(d)$, $\Delta=\Delta_d$, $s\mu+1$ instead of $\mu$, 
$m+i\Delta=z$, we get the expected inclusion
$\Delta_d(F) \supset R_{s\mu+1}(d)$.
\findem

\section{Conclusion of the proofs}

\subsection{Proof of theorems \ref{thr:large} and \ref{thr:strict}}
\label{sec: conclusion de la demo}
Let us denote the stratum $C(E_1,\dots,E_1,\dots,E_r,\dots,E_r)$ by 
$C(E_1^{n_1},\dots,E_r^{n_r})$ where $n_i$ is the number of
copies of $E_i$. According to proposition
\ref{reduction_to_the_computation_of_closure},  
to conclude the proof
of theorem \ref{thr:large} (resp. of theorem \ref{thr:strict}), 
we must show 
that, for $s\geq 2$, $d\geq 2$  and $(s,d)\notin\{(2,2),(2,3),(3,2)\}$
(resp. for $s\geq 1$, $d\geq 2$)
$\overline {C(R_\mu^{s^d})}\supset C(E)$ for some 
staircase $E$ containing $R_{s\mu+1}$ 
(resp. containing $R_{s\mu}$).
By proposition \ref{corollaire: incidence par homotetie}, 
$$\overline {C(R_\mu^{s^d})}\supset C(s.R_{\mu}).$$
Since $s.R_{\mu} \supset R_{s\mu}$, this concludes the proof 
of theorem \ref{thr:strict}. As for theorem \ref{thr:large}, 
taking for $E$ the staircase 
$\Delta_d(\Delta_{d-1}(\dots\Delta_1(s.R_{\mu})))$ 
constructed in lemma \ref{lm:lemme cle}, we have 
$$
\overline {C(s.R_{\mu})} \supset C(E)
$$ 
by proposition 
\ref{prop:incidence pour escalier unique}.
The required inclusion  $\overline {C(R_\mu^{s^d})}\supset
C(E)$ follows immediatly from the last two displayed inclusions. 
\findem

\begin{rem}
\label{rem:induction on dimension}
  The above proof relies heavily on lemma \ref{lm:lemme cle},
which is the key point. This 
key lemma is proved by an induction on the dimension.
\end{rem}

\subsection{The case $r=8,d=3$}
\label{sec:case-r=2d}
The goal of this section is to compute the postulation of 
$8$ fat points of multiplicity $\mu$ in $\PP^3$,
stated in remark \ref{remarque sur le cas 2^d}: 
\begin{prop} \label{prop: le cas r=2^d}
  Let $r=8, d=3$ and $v(d,\delta,\mu^r)=\binom
  {\delta+d}{d}-r \binom{d+\mu-1}{d}$. Then:
  $l(d,\delta,\mu^r)=max(0,v(d,\delta,\mu^r))$. 
\end{prop}
\demo: if $Z_G$ is the generic union of $8$ fat points of 
multiplicity $\mu$, the vector space 
$H^0(I_{Z_G}(\delta))$ being the kernel of the 
restriction morphism: 
\begin{displaymath}
H^0({\cal O}_{\PP^d}(\delta))\fd H^0({\cal
  O}_{Z_G}(\delta)),
\end{displaymath}
its dimension $l(d,\delta,\mu^r)$
is at least 
$$
v(d,\delta,\mu^r)=h^0({\cal O}_{\PP^d}(\delta))-h^0({\cal
  O}_{Z_G}(\delta)).
$$
Let $\Delta=(1,-1,-1)$ and $E=\Delta(2R_{\mu})$.
\\
To prove the reverse inequality 
$l(d,\delta,\mu^r)\leq max(0,v(d,\delta,\mu^r))$,
since $\overline {C(R_\mu^{2^d})} \supset~ C(E)$ 
by proposition \ref{corollaire: incidence par homotetie}
and proposition \ref{prop:incidence pour escalier unique},
it suffices by semi-continuity to exhibit 
a subscheme $Z$ in $C(E)$ such that 
\begin{displaymath}
h^0(I_Z(\delta))=max(0,v(d,\delta,\mu^r)) 
\end{displaymath}
for all $\delta$. Let $\AA^d=\s k[x_1,\dots,x_d] \inc \PP^d$ be an affine
space and $Z$ be the subscheme of $\AA^d$ whose ideal is $I^{E}$. 
By deshomogeneisation, the vector space  $H^0({\cal
  O}_{\PP^d}(\delta))$ is in bijection with the subspace $S_{\delta}\inc 
k[x_1,\dots,x_d]$ containing the polynomials of degree at most
$\delta$, and $H^0(I_Z(\delta))$ corresponds to $I^{E}\cap S_{\delta}$.
Now, $\dim I^{E} \cap S_{\delta}$ is the number of 
monomials in $R_{\delta+1}$ which are not in $E $. Since 
\begin{displaymath}
R_{2\mu} \inc E \inc R_{2\mu+1},
\end{displaymath}
this number is $0$ if
$\delta\leq 2\mu-1$
and $h^0({\cal O}_{\PP^d}(\delta))-\#E$ if $\delta \geq2\mu$..
\findem

\end{document}